\documentclass[letterpaper, paper,11pt]{AAS}		

\usepackage{bm}
\usepackage{amsmath}
\usepackage{subfigure}
\usepackage[colorlinks=true, pdfstartview=FitV, linkcolor=black, citecolor= black, urlcolor= black]{hyperref}
\usepackage{overcite}
\usepackage{footnpag}	
\usepackage{float}

\PaperNumber{26-942}

\allowdisplaybreaks

\begin{document}

\title{Refining GLSDC for Orbit Determination}

\author{Abhijeet\thanks{Graduate Student, Department of Aerospace Engineering, Texas A\& M University, College Station, TX 77843 USA}\footnotemark[2],  
Raja Abhishek Appana\thanks{Graduate Research Assistant, Department of Mechanical and Aerospace Engineering, The Ohio State University, Columbus, OH 43210 USA}\thanks{These authors contributed equally to this work.},
Mrinal Kumar\thanks{Professor and  Elizabeth M. Tinkham Endowed Chair of Aeronautics and Astronautics, Department of Mechanical and Aerospace Engineering. The Ohio State University, Columbus, OH 43210 USA},
\ and Suman Chakravorty\thanks{Professor, Department of Aerospace Engineering, Texas A\& M University, College Station, TX 77843 USA}
}

\maketitle{} 		

\begin{abstract}
This paper investigates the ambiguity of noisy short-arc angles-only orbit determination for a Molniya-orbit case study. Conventional Gauss Initial Orbit Determination followed by Gaussian Least-Squares Differential Correction is first shown to be unreliable under high-noise short-arc conditions, frequently failing to converge or converging to hyperbolic local solutions. To explore the candidate solution space, a Lambert-based initialization procedure is used over a grid of assumed ranges and angular-observation pairs, and the resulting states are refined using nonlinear least-squares optimization. The unconstrained solution set reveals multiple orbit families, including reentry elliptic, bounded elliptic, xGEO elliptic, and hyperbolic trajectories, all of which can reproduce the observed angular arc with comparable residuals. The solutions exhibit a structured range–velocity relationship, indicating that the angular measurements primarily constrain apparent line-of-sight motion rather than absolute range. Finally, physically motivated constraints are introduced to isolate Molniya-like bounded elliptic solutions. The results demonstrate that short-arc angles-only orbit determination is fundamentally non-unique without a priori information and that constraints should be interpreted as orbit-family selection tools rather than proof of uniqueness.
\end{abstract}

\section{Introduction}

Accurate orbit determination is critical for space situational awareness, particularly given the recent surge in resident space objects (RSOs) and their associated collision risks. Although the global sensor network has expanded, limitations in field of view and external disturbances result in data being available only over too-short arcs (TSAs). Small-scale optical telescopes often introduce high measurement noise ($> 15\text{ arcsec}$) and can only capture partial state measurements, further complicating definitive orbit determination.

Extensive research has focused on defining Admissible Regions (ARs) from TSA measurements. First proposed by Milani et al. \cite{milani2004orbit} for tracking very short-arc celestial objects, the AR concept defines the complete set of orbits that satisfy available observations within an acceptable error bound. This framework has since been widely applied to asteroids \cite{milani2005astrometry}, space debris \cite{tommei2007orbit, popplewell2025tracking}, and active satellites. It has been extended to include Constrained Admissible Regions (CAR) \cite{kelecy2013application} and Probabilistic Admissible Regions (PAR) \cite{hussein2014probabilistic, mishra2024geometric} as well. The AR is typically mapped onto the range and range-rate $(\rho, \dot{\rho})$ plane and bounds the state space based on physical or geometric criteria. Various techniques, including multi-hypothesis filtering \cite{stauch2018joint} and ranging methods \cite{tao2020uncertainty}, are then deployed to identify the true orbit.

To determine the admissible region, an attributable vector consisting of the angles and angle rates at a point of measurement is made available. This vector is then used to determine the curve of zero energy ($\varepsilon = 0$) on the range-range-rate plane \cite{demars2012initial}. Assuming an earth-orbiting object, this curve defines the admissible region inside which all orbits with negative energy exist ($\varepsilon < 0$). If a priori orbital characteristic information is available, then the region can be further restricted to CAR \cite{demars2013probabilistic}. It needs to be noted that the angle rate information is rarely available from small-scale sensors, which typically provide only time-series angular data. Hussein et al. \cite{hussein2018probabilistic} derive the attributable vector from these angles-only data using the method proposed by Maruskin et al. \cite{maruskin2009correlation}. The work proposes a quadratic kinematic model of the angles, and the least-squares approach is employed to approximate the angular rate and acceleration, which inherently introduces systematic errors in the data.

Furthermore, the vast majority of existing literature in this area assumes either noise-free data or data with a noise level (standard deviation) under $10\text{ arcsec}$  \cite{popplewell2025tracking, hussein2018probabilistic}. This creates a significant capability gap as small-scale telescopes have angular standard deviation noise exceeding $15\text{ arcsec}$.

In this paper, we first demonstrate the limitations of conventional initial orbit determination (IOD) and subsequent orbit determination (OD) techniques when applied to noisy short-arc angles-only data. Under such conditions, the observation geometry provides only limited information about the range and range-rate, and classical methods may either fail to converge or converge to one of several possible local solutions. We therefore investigate the structure of the candidate solution space rather than assuming that a unique orbit can be recovered from the available measurements.

For a short observation arc of an elliptic reference orbit, initial guesses are generated by solving Lambert's problem over varying range assumptions and different pairs of angular observations. These initial guesses are then refined through nonlinear optimization. The resulting candidate solutions are not restricted to the original orbit family; instead, they may correspond to reentry elliptic orbits, bounded elliptic orbits between the reentry and xGEO regimes, xGEO elliptic orbits, or even hyperbolic trajectories. In the absence of additional a priori information, all of these candidate orbits can represent the observed angles over the short arc with comparable residuals.

The analysis shows that these solutions exhibit a structured range--velocity relationship: as the magnitude of the estimated position increases, the corresponding velocity magnitude also increases. This behavior arises because the angles-only measurements primarily constrain the apparent angular motion, or angular-rate pattern, rather than the absolute range. Consequently, distinct orbit families can reproduce similar line-of-sight histories over the short observation interval. Finally, we show that when physically motivated constraints or a priori information are introduced, the optimization can be guided toward a solution close to the original Molniya orbit. Thus, the proposed analysis highlights both the ambiguity inherent in short-arc angles-only orbit determination and the role of a priori information in selecting a physically meaningful solution.

The remainder of the paper is organized as follows. The ``Background" section provides the necessary background on GLSDC and Lambert's problem. Following which, the ``Methodology" section describes the simulation setup, including the orbital parameters, measurement data, noise model, and orbit-determination procedure. This section also briefly examines the use of angular rates to obtain an approximate range estimate. We provide the results of the proposed approach in the ``Results" section. It includes the classification of the recovered orbit families and the use of a priori range information to improve convergence toward the reference orbit. Finally, the concluding section summarizes the main findings and outlines possible directions for future work.

\section{Background}\label{sec:background}

This section gives a brief overview of the established IOD and OD methods that are used in this paper. To show the limitations of classical methods, the angles-only Gauss IOD method is used to determine the initial orbit, which is further refined using Gaussian Least-Squares Differential Correction (GLSDC). For generating initial guesses to be used in the constrained optimization setup, the problem is posed as a Lambert's problem with varying range values. An outline of each method is provided in the rest of this section.

\subsection{Gauss's method for Initial Orbit Determination \cite{curtis2019orbital}}

Gauss’s method requires three time-ordered angular observations (line-of-sight unit vectors $\hat{\mathbf{L}}_1, \hat{\mathbf{L}}_2$ $, \hat{\mathbf{L}}_3$) at epochs $t_1, t_2, t_3$. The observer’s topocentric position vectors $\mathbf{R}_i$ are assumed to be known. The fundamental relation connects the geocentric position vector $\mathbf{r}_i$ to the unknown range $\rho_i$ via:
\begin{equation} \label{eq: gauss_range}
    \mathbf{r}_i = \mathbf{R}_i + \rho_i \hat{\mathbf{L}}_i \quad (i = 1, 2, 3)
\end{equation}

The core insight behind Gauss's method is that the orbital motion is coplanar, meaning that the three position vectors are linearly dependent. This is mathematically expressed as:
\begin{equation} \label{eq: gauss_coplanar}
    c_1 \mathbf{r}_1 - c_2 \mathbf{r}_2 + c_3 \mathbf{r}_3 = \mathbf{0}
\end{equation}

where the coefficients $c_1, c_2, c_3$ are functions of the time intervals $\tau$ and the changing orbital geometry. These time intervals are typically defined as:$$\tau_1 = t_1 - t_2, \quad \tau_3 = t_3 - t_2, \quad \tau_2 = t_3 - t_1$$

By approximating these coefficients using a Taylor series expansion in terms of the gravitational parameter $\mu$ and the central radius $r_2 = \|\mathbf{r}_2\|$, Gauss expressed them as:
\begin{align}
    c_1 &\approx \frac{\tau_3}{\tau_2} \left(1 + \frac{\mu}{6r_2^3}(\tau_2^2 - \tau_3^2)\right) \\
    c_3 &\approx \frac{\tau_1}{\tau_2} \left(1 + \frac{\mu}{6r_2^3}(\tau_2^2 - \tau_1^2)\right) \\
    c_2 &= 1
\end{align}

By substituting Eq. \ref{eq: gauss_range} into Eq. \ref{eq: gauss_coplanar} and taking the dot and cross products with the line-of-sight vectors, the problem is reduced to solving for the middle range $\rho_2$ and the middle radius $r_2$. This yields a system of two coupled equations. The first is a linear relationship derived from the geometry:$$\rho_2 = A + \frac{B}{r_2^3}$$(where $A$ and $B$ are scalar constants computed purely from the known observer vectors $\mathbf{R}_i$, line-of-sight vectors $\hat{\mathbf{L}}_i$, and time intervals $\tau_i$). The second relation is the fundamental geometric identity:$$r_2^2 = \|\mathbf{R}_2 + \rho_2 \hat{\mathbf{L}}_2\|^2 = \rho_2^2 + 2\rho_2 (\mathbf{R}_2 \cdot \hat{\mathbf{L}}_2) + \|\mathbf{R}_2\|^2$$Substituting the linear $\rho_2$ expression into the geometric identity results in the Gauss eighth-order polynomial in terms of $r_2$:
\begin{equation} \label{eq: gauss_iod}
    r_2^8 + a r_2^6 + b r_2^3 + c = 0
\end{equation}

\subsection{Gaussian Least-Squares Differential Correction \cite{schutz2004statistical}}

The Gaussian Least Squares Differential Correction (GLSDC) is a recursive process to estimate an object's orbital state (position and velocity) from a series of measurements. It iteratively refines an initial guess until a convergence criterion is met. The process consists of the following steps, which are repeated until convergence:
\begin{enumerate}
    \item \textbf{Initial Guess (IOD):}
    Obtain an initial estimate of the state vector, $\hat{\mathbf{x}}_0^{(0)} = [\hat{\mathbf{r}}_0, \hat{\mathbf{v}}_0]^T$, at an initial time $t_0$. This is typically done using an Initial Orbit Determination (IOD) method, such as Gauss's method for angles-only data.

    \item \textbf{Propagate State and STM:}
    Numerically integrate the equations of motion and the variational equations from the current state estimate $\hat{\mathbf{x}}_0$ to each measurement time $t_i$:
    \begin{align*}
        \dot{\mathbf{x}}(t) &= \mathbf{f}(\mathbf{x}(t)) = 
        \begin{bmatrix}
            \mathbf{v}(t) \\
            -\frac{\mu}{r^3}\mathbf{r}(t)
        \end{bmatrix} \\
        \dot{\Phi}(t, t_0) &= F(t)\Phi(t, t_0), \quad \text{with } \Phi(t_0, t_0) = I_{6 \times 6}
    \end{align*}
    where $F(t) = \frac{\partial \mathbf{f}}{\partial \mathbf{x}}$ is the Jacobian of the system dynamics.

    \item \textbf{Compute Predicted Measurements:}
    Using the propagated states $\hat{\mathbf{x}}(t_i)$, calculate the predicted measurements $\hat{\mathbf{y}}_i$ (e.g., azimuth and elevation) for each observation time using the observation model $\mathbf{h}(\cdot)$:
    \[ \hat{\mathbf{y}}_i = \mathbf{h}(\hat{\mathbf{x}}(t_i)) \]

    \item \textbf{Compute Residuals:}
    Calculate the observation residual vector $\Delta \mathbf{y}_c$ by subtracting the estimated measurements from the actual measurements $\tilde{\mathbf{y}}$:
    \[ \Delta \mathbf{y}_c = \tilde{\mathbf{y}} - \hat{\mathbf{y}} \]

    \item \textbf{Compute the Sensitivity Matrix (H):}
    The matrix $H$ relates changes in the initial state to changes in the observations. It is constructed by stacking the sensitivity for each measurement time $t_i$:
    \[ H = 
    \begin{bmatrix}
        \frac{\partial \mathbf{h}_1}{\partial \mathbf{x}_1} \Phi(t_1, t_0) \\
        \vdots \\
        \frac{\partial \mathbf{h}_m}{\partial \mathbf{x}_m} \Phi(t_m, t_0)
    \end{bmatrix}
    \]
    For angles-only measurements (azimuth $A$ and elevation $a$), the matrix $\frac{\partial \mathbf{h}_i}{\partial \mathbf{x}_i}$ at a given time has the form:
    \[
    \frac{\partial \mathbf{h}}{\partial \mathbf{x}} = 
    \begin{bmatrix}
        \frac{\partial A}{\partial r_x} & \frac{\partial A}{\partial r_y} & \frac{\partial A}{\partial r_z} & 0 & 0 & 0 \\
        \frac{\partial a}{\partial r_x} & \frac{\partial a}{\partial r_y} & \frac{\partial a}{\partial r_z} & 0 & 0 & 0
    \end{bmatrix}
    \]
    The full $H$ matrix for $m$ measurements is a $2m \times 6$ matrix.

    \item \textbf{Calculate the Correction:}
    Solve the linear normal equation for the correction to the initial state, $\Delta \mathbf{x}_{0c}$:
    \[ \Delta \mathbf{x}_{0c} = (H^T W H)^{-1} H^T W \Delta \mathbf{y}_c \]
    where $W$ is a weighting matrix, often related to the inverse of the measurement noise covariance.

    \item \textbf{Update State and Check Convergence:}
    Update the initial state estimate:
    \[ \hat{\mathbf{x}}_0^{(k+1)} = \hat{\mathbf{x}}_0^{(k)} + \Delta \mathbf{x}_{0c} \]
    The process is terminated when a convergence criterion is met, typically by checking for a minimal change in the cost function $J_c = \frac{1}{2} \Delta \mathbf{y}_c^T W \Delta \mathbf{y}_c$. If not converged, return to Step 2.
\end{enumerate}

\subsection{Lambert's problem \cite{bate2020fundamentals}}
Lambert’s problem is formulated as a two-point boundary-value problem with known time-of-flight in Keplerian mechanics. Given two time-tagged geocentric position vectors, $\mathbf{r}_1$ at $t_1$ and $\mathbf{r}_2$ at $t_2$, the objective is to determine the unique Keplerian trajectory connecting them over the specified time-of-flight (TOF):$$\Delta t = t_2 - t_1$$
The fundamental physics is governed by the two-body differential equation:$$\ddot{\mathbf{r}} + \frac{\mu}{\|\mathbf{r}\|^3}\mathbf{r} = \mathbf{0}$$
For our current application, to handle elliptic ($a > 0$), parabolic ($a \to \infty$), and hyperbolic ($a < 0$) trajectories within a single mathematical framework, the universal Lambert's solution is used to determine the initial guess of the orbit for a given range value. For a comprehensive understanding of the universal solution, the readers are encouraged to refer to \textit{Fundamentals of Astrodynamics}\cite{bate2020fundamentals}.

\section{Methodology} \label{sec:methodology}
This section describes the methodology used to investigate the short-arc angles-only orbit determination problem for a Molniya orbit. The objective is not only to estimate the orbital state from limited angular measurements, but also to understand the structure of the resulting optimization problem. In particular, we first examine the performance of the classical Gaussian Least-Squares Differential Correction (GLSDC) method for short observation arcs and show that it may either fail to converge or converge to an incorrect local solution. This behavior motivates a deeper analysis of the admissible solution space, where multiple candidate orbits can produce similar angular residuals over a short time interval.

To study this ambiguity, we formulate the orbit determination problem as a constrained nonlinear minimization problem. The measurement model is based on angles-only observations, and the numerical experiments are performed using a prescribed Molniya orbit with specified orbital elements, observation geometry, sampling interval, and measurement noise. The resulting objective function is minimized, subject to physically meaningful orbital constraints, using the interior-point algorithm in MATLAB’s \texttt{fmincon}. Since the short-arc problem is sensitive to the initial estimate, multiple initial conditions are generated by solving Lambert’s problem over the observation interval. These Lambert-based initial guesses allow the optimizer to explore different regions of the solution space and reveal the existence of multiple local optima. The final candidate solutions are then analyzed and classified according to their orbital parameters, residual cost, and physical feasibility.

\subsection{Simulation and Observation Setup}
As noted previously, the orbit considered in this study is a Molniya orbit. The corresponding orbital parameters, together with the true anomaly at the initial observation epoch ($t = 0$), are summarized in Table \ref{tab:orbit_parameters}. The ground-based observatory is assumed to be located in Hawaii, and the relevant site parameters, along with the observation timing, are provided in Table \ref{tab:loc_obs}. The measurements consist of azimuth and elevation angles sampled approximately every $10$ seconds. A total of $20$ measurements are collected, resulting in an observation arc of approximately $200$ seconds. The measurement configuration is summarized in Table \ref{tab: measurement}. The synthetic measurement data are generated from the reference Molniya orbit, after which zero-mean noise with a standard deviation of $18 \text{ arcseconds}$ is added to obtain a more realistic observation set. Therefore, the measurements reported in Table \ref{tab: measurement} correspond to the noisy azimuth and elevation data used in the orbit determination process.

\begin{table}[!htbp]
    \centering
    \begin{tabular}{|c|c|c|c|c|c|c|}
    \hline
       \textbf{Orbital Parameters} & a  & e & $i$ & $\Omega$ & $\omega$ & $\nu$ \\
       \hline
       \textbf{Value} & 2.6559e+07  & 0.6757 & $63.5^{o}$ & $167.25^{o}$ & $277.98^{o}$ & $240^{o}$\\
       \hline
    \end{tabular}
    \caption{Orbital Parameters of the Molniya Orbit and the true anomaly at the first point of observation}
    \label{tab:orbit_parameters}
\end{table}

\begin{table}[!htbp]
    \centering
    \begin{tabular}{|c|c|c|c|c|}
        \hline
        Latitude & Longitude & Altitude & Local Date & Local Time\\
        \hline
        $20^{o} 42' 30''$ & $-156^{o} 15' 29''$ & 3052 m & 10-1-2025 & 2:00\\
        \hline
    \end{tabular}
    \caption{Location of the observation site and the time of observation}
    \label{tab:loc_obs}
\end{table}

\begin{table}[h]
    \centering
    \begin{tabular}{|c|c|c|}
    \hline
       t (seconds)  & Azimuth ($A$ (degrees))
         & Elevation ($a$ (degrees))\\
         \hline
         0 & -80.384 & 20.696\\
         \hline
         10.774 & -80.499 & 20.668\\ 
         \hline
         21.548 & -80.606 & 20.615
       \\
       \hline
      32.323 & -80.722 & 20.561\\
      \hline
      43.097 & -80.829 & 20.524\\
      \hline
      53.871 & -80.931 & 20.479\\
      \hline
      64.645 & -81.038 & 20.439\\
      \hline
       75.42 & -81.15 & 20.392\\ 
       \hline
      86.194 & -81.266 & 20.345\\
      \hline
      96.968 & -81.364 & 20.312\\ 
      \hline
      107.74 & -81.483 & 20.253\\
      \hline
      118.52 & -81.595 & 20.22\\
      \hline
      129.29 & -81.694 & 20.161\\
      \hline
      140.06 & -81.801 & 20.129\\
      \hline
      150.84 & -81.917 & 20.079\\
      \hline
      161.61 & -82.031 & 20.038\\
      \hline
      172.39 & -82.157 & 19.988\\
      \hline
      183.16 & -82.258 & 19.936\\
      \hline
      193.94 & -82.361 & 19.897\\
      \hline
      204.71 & -82.479 & 19.851\\
      \hline
    \end{tabular}
    \caption{Azimuth and elevation data considered for orbit determination}
    \label{tab: measurement}
\end{table}

\begin{figure}[htbp!]
    \centering
    \includegraphics[width=0.7\linewidth]{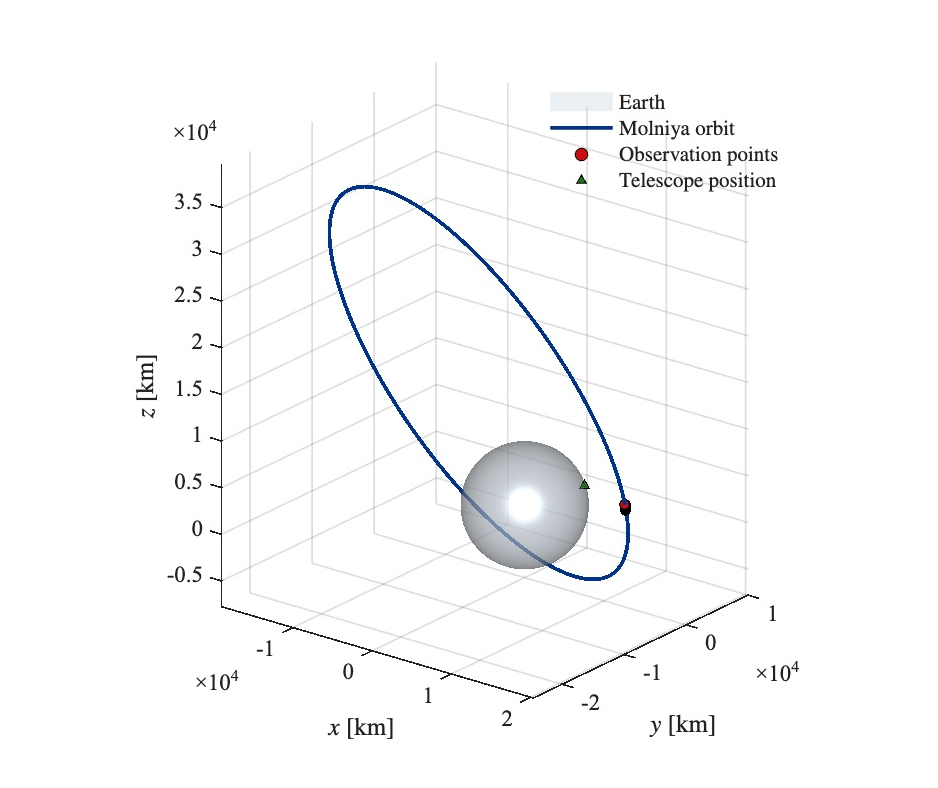}
    \caption{The Molniya orbit along with the telescope position and the region of observation}
    \label{fig:molniya}
\end{figure}

\subsection{Initial Orbit Determination}
To generate a diverse set of initial guesses for the orbit determination problem, we employ a Lambert-based initialization strategy. Two angular measurements are first selected, consisting of azimuth and elevation pairs $((A_{t_1}, a_{t_1}))~\text{and}~((A_{t_2}, a_{t_2}))$. For each selected pair, a candidate range is assigned, allowing the corresponding position vectors to be constructed along the two line-of-sight directions. A Lambert problem is then solved between these two position vectors to obtain an initial state estimate $([\mathbf{r}_{t_1}, \mathbf{v}_{t_1}])$, which is subsequently used as the initial guess for the orbit determination problem.

In this study, the first measurement is paired with each of the remaining 19 measurements, resulting in 19 possible angular-measurement combinations for initial orbit determination. In addition, the assumed range is varied from 100 km to 1000 km in increments of 100 km, and from 1000 km to 30000 km in increments of 1000 km, avoiding duplicate range values at 1000 km. This produces a total of 39 candidate range values. Therefore, the initialization procedure generates $(39 \times 19 = 741)$ Lambert-based initial orbit estimates. However, one should note that most of these initial guesses will be very poor as the range at the time of observation is approximately 21000 Kms.

\subsection{Orbit Determination}

The orbit determination problem is formulated as a nonlinear least-squares optimization problem in which the initial Cartesian state,
\[
\mathbf{x}_0 =
\begin{bmatrix}
\mathbf{r}_0^T & \mathbf{v}_0^T
\end{bmatrix}^T,
\]
is estimated from a set of angles-only measurements. The measurement vector consists of the observed azimuth and elevation angles collected over the short observation arc. For a given initial state ($\mathbf{x}_0$), the orbit is propagated to each observation time, and the corresponding predicted azimuth and elevation angles are computed from the relative geometry between the object and the ground-based observatory. The resulting objective function is written as

\[
J(\mathbf{x}_0)=\frac{1}{2}
\left(
\tilde{\mathbf{y}} - \mathbf{f}(\mathbf{x}_0)
\right)^T
W
\left(
\tilde{\mathbf{y}} - \mathbf{f}(\mathbf{x}_0)
\right),
\]
where ($\tilde{\mathbf{y}}$) is the vector of measured azimuth and elevation angles, $(\mathbf{f}(\mathbf{x}_0)) $is the nonlinear mapping from the initial state to the predicted measurements, and (W) is the weighting matrix.

The classical GLSDC method solves this problem by linearizing the measurement model about the current estimate and iteratively applying a correction based on the measurement sensitivity matrix. However, for short-arc angles-only measurements, the resulting optimization landscape can be highly nonconvex and weakly constrained in range and range-rate. Consequently, GLSDC may fail to converge or may converge to an incorrect local solution that produces a small angular residual but corresponds to an inaccurate orbit.

To investigate this behavior, the same orbit determination problem is also solved using the interior-point algorithm in MATLAB's \texttt{fmincon}. Unlike the standard GLSDC update, the interior-point framework allows the optimization problem to be solved both with and without explicit orbital constraints. In addition, the Hessian approximation used in the Newton-type formulation includes the second-order correction term associated with the nonlinear measurement model. Specifically, if
\[
H = \frac{\partial \mathbf{f}}{\partial \mathbf{x}}
\]
is the first-order measurement sensitivity matrix and
\[
G = \frac{\partial^2 \mathbf{f}}{\partial \mathbf{x}_i \partial \mathbf{x}_j}
\]
denotes the second derivative tensor of the measurement model, the Newton-type Hessian contribution can be written as
\[
\nabla^2 J
\approx
H^T W H
-G \otimes W
\left(
\tilde{\mathbf{y}} - \mathbf{f}(\mathbf{x}_0)
\right).
\]
This additional term accounts for the local curvature of the nonlinear measurement function and is consistent with the second-order expansion of the least-squares cost. The interior-point optimizer can use this Hessian information to improve the local search direction when solving the nonlinear optimization problem.

\begin{figure}[h]
    \centering
    \includegraphics[width=0.7\linewidth]{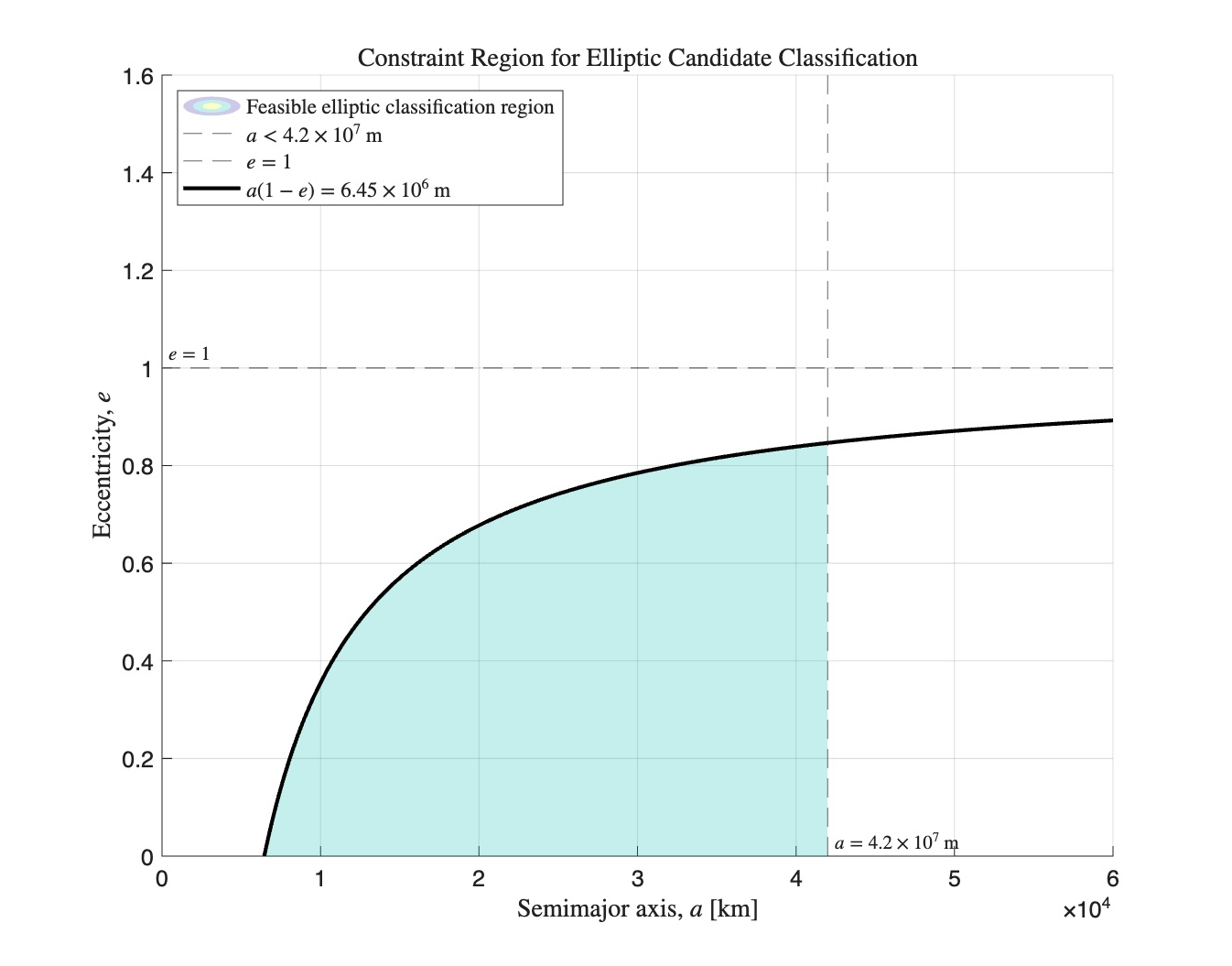}
    \caption{A depiction of applied constraints in the semimajor axis- eccentricity plane}
    \label{fig:constraints}
\end{figure}

Two optimization cases are considered. In the first case, the orbit determination problem is solved without imposing orbital constraints. This unconstrained formulation is used to expose the full structure of the solution space and to identify candidate optima that may correspond to elliptic, hyperbolic, reentry, or high-altitude solutions. In the second case, a constrained optimization problem is solved to target physically plausible elliptic solutions. The constrained problem is given by
\[
\min_{\mathbf{x}_0} J(\mathbf{x}_0)
\]
subject to
\[
\epsilon < 0,
\]
\[
a < 4.2 \times 10^7 \ \text{m},
\]
and
\[
a(1-e) > 6.45 \times 10^6 \ \text{m},
\]
where, $\epsilon$ is the specific orbital energy, $a$ is the semimajor axis, and $e$ is the eccentricity. The first constraint enforces a bound elliptic orbit, the second constraint restricts the semimajor axis to remain below the approximate geosynchronous radius, and the third constraint ensures that the periapsis radius remains above the Earth reentry/impact region. These constraints are depicted in the semimajor-eccentricity plane in Figure \ref{fig:constraints}.

It is important to note that these constraints are not used to discard the other solution families from the analysis. Hyperbolic solutions, reentry solutions, and solutions with semimajor axes beyond the geosynchronous radius are still retained in the results obtained from the unconstrained optimization and are included in the classification of local optima. The constrained formulation is instead used to isolate the region of the solution space in which physically plausible elliptic solutions are expected to have a higher likelihood. Therefore, the constraints serve both as an optimization tool and as a classification criterion for distinguishing Molniya-compatible candidate solutions from other mathematically admissible short-arc fits.

Because the short-arc angles-only problem is highly sensitive to the initial condition, the optimization is initialized using the Lambert-based procedure described previously. Each Lambert solution provides a candidate initial state $([\mathbf{r}_{t_1},\mathbf{v}_{t_1}])$, which is then passed to \texttt{fmincon}. By solving the constrained and unconstrained optimization problems over a large set of Lambert-generated initial guesses, the method explores multiple basins of attraction in the objective function and enables classification of the resulting local optima.

\section{Results}\label{sec:results}

This section presents the results obtained from the short-arc angles-only orbit determination problem. The goal is to examine the structure of the solution space obtained from unconstrained optimization and to compare it with the solution set obtained after imposing physically motivated constraints. The results show that the short observation arc admits multiple orbit families with comparable angular residuals. Therefore, the solution of the orbit determination problem should not be interpreted as unique based only on the final cost. Before delving into that, we present an overview of solutions obtained from GLSDC.

\subsection{Gaussian Least-Squares Differential Correction}
In this subsection, we briefly discuss the results obtained using the GLSDC approach and explain why the subsequent analysis is carried out using \texttt{fmincon}. In particular, we show that while GLSDC provides useful baseline estimates, it is not well suited for systematically exploring the multiple solution families present in the short-arc angles-only problem. This motivates the use of \texttt{fmincon}, which allows the optimization problem to be solved with explicit constraints and from a large set of initial guesses.

To begin the analysis, we first applied the classical Gauss IOD procedure followed by GLSDC refinement. The measurement set consists of 20 azimuth--elevation pairs, and Gauss IOD requires three angular measurements to generate an initial orbit estimate. Therefore, all possible $( {20 \choose 3}=1140 )$ measurement triplets were used to initialize the GLSDC algorithm.

Out of these 1140 initial guesses, 472 cases failed to converge. The dominant causes of failure were ill-conditioned or non-invertible Hessian matrices and excessively large update steps, which caused the cost to increase during the iterative process. These cases were therefore discarded from the subsequent analysis. The remaining 668 cases converged; however, all of the converged solutions corresponded to hyperbolic orbits. The observed range of semimajor-axis values for these hyperbolic solutions was from $(-4.3793\times 10^{-18})$ m to $ (-3.9744\times 10^{7})$ m.

\begin{figure}
    \centering
    \includegraphics[width=0.8\linewidth]{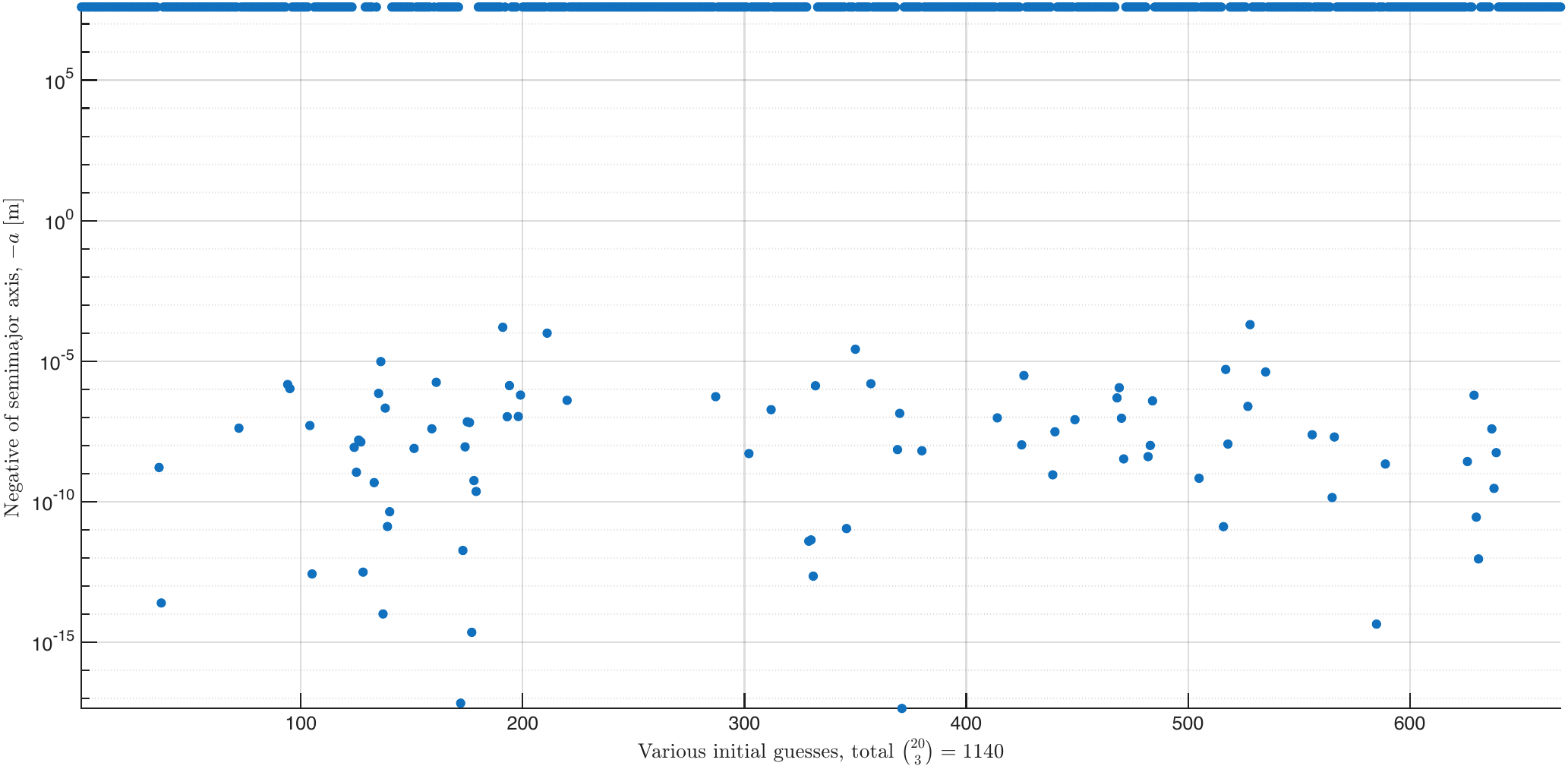}
    \caption{Distribution of hyperbolic candidate solutions obtained from different Lambert-based initial guesses. The negative semimajor axis, \(-a\), is shown on a logarithmic scale for all 668 initial guesses which converged, highlighting the spread of open-orbit solutions produced by the short-arc angles-only measurements.}
    \label{fig:GLSDC_fail}
\end{figure}

Figure~\ref{fig:GLSDC_fail} shows the variation of the negative semimajor axis for the 668 converged GLSDC solutions. The horizontal axis represents the different converged initial guesses, while the vertical axis shows $(-a)$ on a logarithmic scale. It can be seen that most of the converged solutions collapse to approximately the same semimajor-axis value, corresponding to $(a=-3.9744\times 10^{7})$ m. The remaining orbital elements were also found to be similar for this dominant solution family. Thus, although multiple Gauss-IOD initial guesses were used, the GLSDC procedure primarily converged to a single hyperbolic solution branch rather than revealing the broader class of possible orbits.

This behavior indicates that GLSDC is not sufficient for systematically exploring the multiple local minima present in the short-arc angles-only orbit determination problem. In particular, the method either fails to converge for a significant fraction of the initial guesses or converges predominantly to a hyperbolic solution. This motivates the use of Lambert-based initialization over varying range assumptions, followed by nonlinear optimization using \texttt{fmincon}. As shown in the following subsections, this approach allows a richer set of candidate solutions to be recovered, including reentry elliptic, bounded elliptic, xGEO elliptic, and hyperbolic orbit families.

\subsection{Unconstrained Optimization Results}

\begin{figure}[htbp]
\centering
\subfigure[All elliptic orbits.]{
    \includegraphics[width=0.48\textwidth]{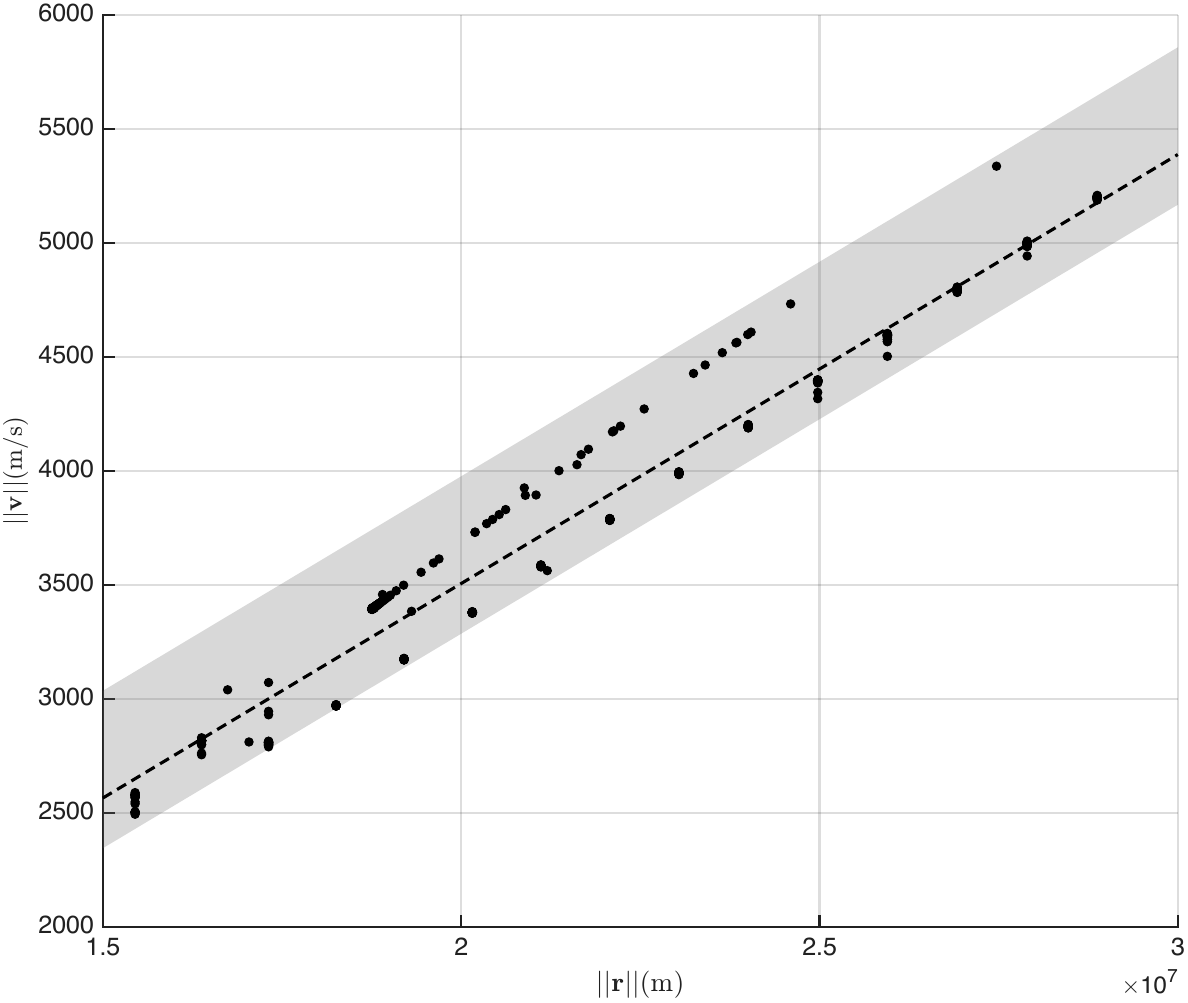}
    \label{fig:elliptic_unc}
}
\hfill
\subfigure[All hyperbolic orbits.]{
    \includegraphics[width=0.48\textwidth]{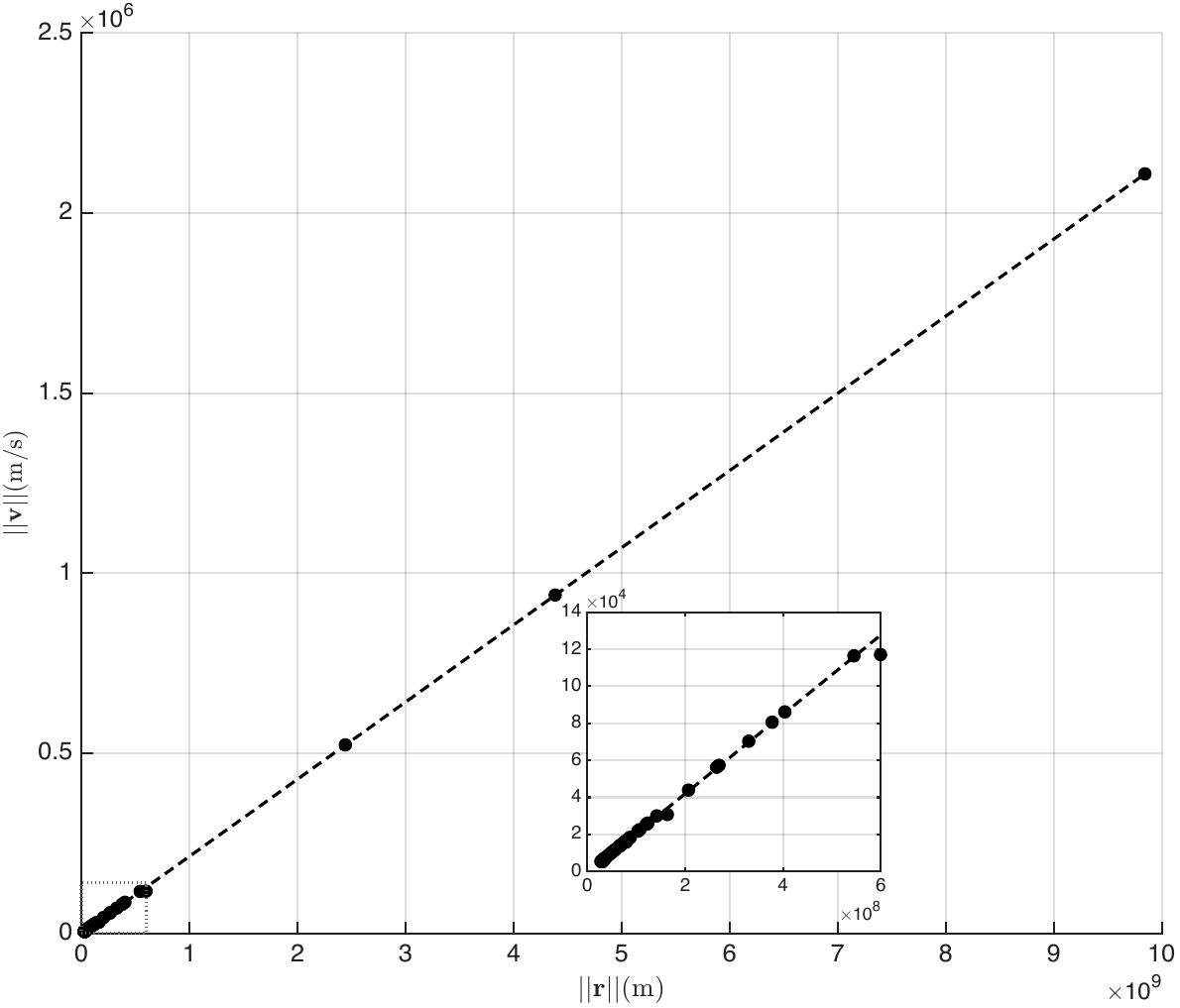}
    \label{fig:hyperbolic_unc}
}

\caption{Distribution of unconstrained orbit-determination solutions in the $||\mathbf{v}||$ -versus- $||\mathbf{r}||$ plane. The elliptic solutions occupy a narrow, approximately linear region, highlighted in gray, while the hyperbolic solutions lie close to a distinct linear trend. This structure indicates that the short-arc angles-only data admit multiple range--velocity combinations with similar angular residuals.}
\label{fig:elliptic_hyperbolic}

\end{figure}

The unconstrained optimization problem was first solved using the Lambert-generated initial guesses described in the previous section. The resulting solutions were classified into closed elliptic orbits and open hyperbolic orbits based on their orbital energy. Figure~\ref{fig:elliptic_hyperbolic} shows the corresponding solutions in the $(||\mathbf{v}||)$-versus-$(||\mathbf{r}||)$ plane.

The elliptic solutions, shown in Figure~\ref{fig:elliptic_unc}, do not lie exactly on a single straight line. However, they are concentrated within a narrow, approximately linear patch, which is highlighted in gray. The hyperbolic solutions, shown in Figure~\ref{fig:hyperbolic_unc}, lie close to a more distinct linear trend. This behavior indicates that the short-arc angles-only data allow multiple combinations of position and velocity magnitude to fit the same observation set. In particular, the angular measurements constrain the line-of-sight geometry more strongly than they constrain the absolute range and range-rate.

\begin{figure}[htbp]
\centering
\includegraphics[width=0.7\linewidth]{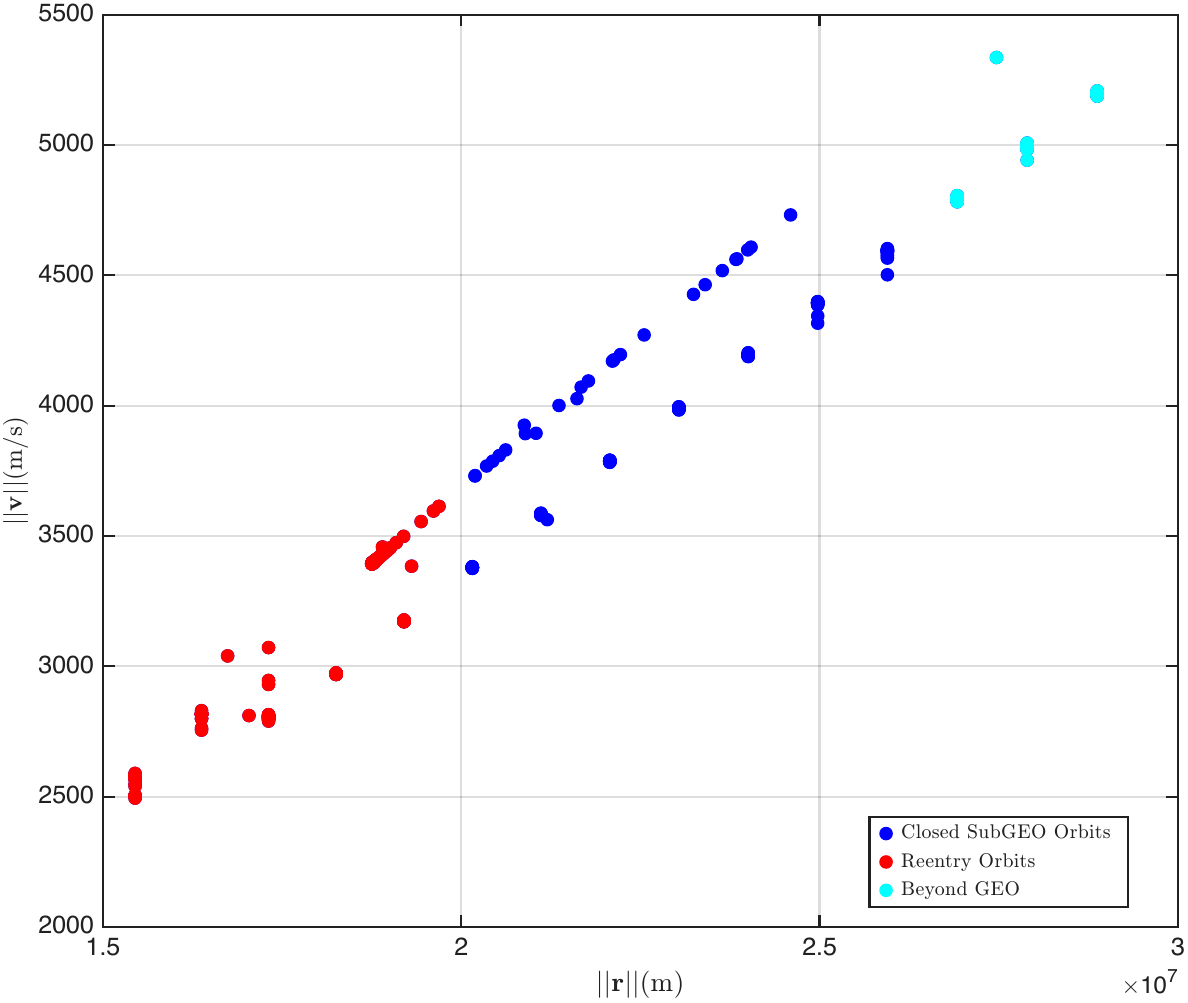}
\caption{Classification of unconstrained elliptic solutions in the $(||\mathbf{v}||)$-versus-$(||\mathbf{r}||)$ plane. Reentry solutions are shown in red, bounded elliptic solutions below the xGEO threshold are shown in blue, and xGEO solutions satisfying $(a > 4.2\times 10^{7}) m$ are shown in cyan. The ordering of the clusters shows that larger estimated ranges are associated with larger velocity magnitudes.}
\label{fig:class_elliptic}
\end{figure}

The elliptic solutions were further classified into three groups, as shown in Figure~\ref{fig:class_elliptic}. Reentry orbits are shown in red, xGEO orbits satisfying $(a > 4.2\times 10^{7})$ m are shown in cyan, and bounded elliptic orbits that are neither reentry nor xGEO are shown in blue. The reentry orbits occupy the lowest range of both $(||\mathbf{r}||)$ and $(||\mathbf{v}||)$. The bounded elliptic orbits appear at larger position and velocity magnitudes, while the xGEO orbits appear at still higher values.

This ordering shows that the unconstrained optimization produces a structured family of solutions. The closer orbits have smaller velocity magnitudes, while larger estimated ranges require larger velocity magnitudes. As the estimated range increases, the solution family transitions from reentry elliptic orbits to bounded elliptic orbits, then to xGEO elliptic orbits, and finally to open hyperbolic orbits.

Table~\ref{tab:range} summarizes the range of orbital elements obtained from the unconstrained solution set. The inclination and the right ascension of the ascending node remain relatively close to the true Molniya orbit values. In contrast, the argument of periapsis and the initial true anomaly show a much larger spread. However, the final column shows that the combined angle $(\omega+\nu_0) $remains close to the true value.

This result suggests that the short-arc measurements constrain the angular location of the object along the orbit more strongly than they constrain the individual values of $(\omega) and (\nu_0)$. Thus, different combinations of argument of periapsis and true anomaly can produce nearly equivalent short-arc line-of-sight histories.

\begin{table}[!htbp]
\centering
\begin{tabular}{|p{2.0cm}|p{2.2cm}|p{1.1cm}|p{1.3cm}|p{1.3cm}|p{1.3cm}|p{1.3cm}|p{1.4cm}|}
\hline
\textbf{Orbital Parameters} & a (m) & e & i & $\Omega$ & $\omega$ & $\nu_{0}$ & $\omega + \nu_{0}$ \\
\hline
\textbf{Actual Value} & $2.6559\times 10^{7}$ & 0.6757 & $63.5^{\circ}$ & $167.25^{\circ}$ & $277.98^{\circ}$ & $240^{\circ}$ & $517.98^{\circ}$\\
\hline
\textbf{Observed min.} & $8.7849\times 10^{6}$ & 0.38369 & $63.13^{\circ}$ & $159.28^{\circ}$ & $206.84^{\circ}$ & $186.74^{\circ}$ & $516.71^{\circ}$\\
\hline
\textbf{Observed max.} & $8.1329\times 10^{8}$ & $0.97302$ & $72.214^{\circ}$ & $174.84^{\circ}$ & $330.31^{\circ}$ & $313.38^{\circ}$ & $520.22^{\circ}$\\
\hline
\end{tabular}
\caption{Reference Molniya orbital parameters and the range of orbital elements recovered from the unconstrained solution set. The inclination and right ascension of the ascending node remain relatively close to the true values, while $(\omega)$ and $(\nu_0)$ vary over a wider range. However, the combined angle $(\omega+\nu_0)$ remains tightly constrained.}
\label{tab:range}
\end{table}

\begin{figure}[htbp!]
\centering
\subfigure[Reentry orbits.]{
    \includegraphics[width=0.48\textwidth]{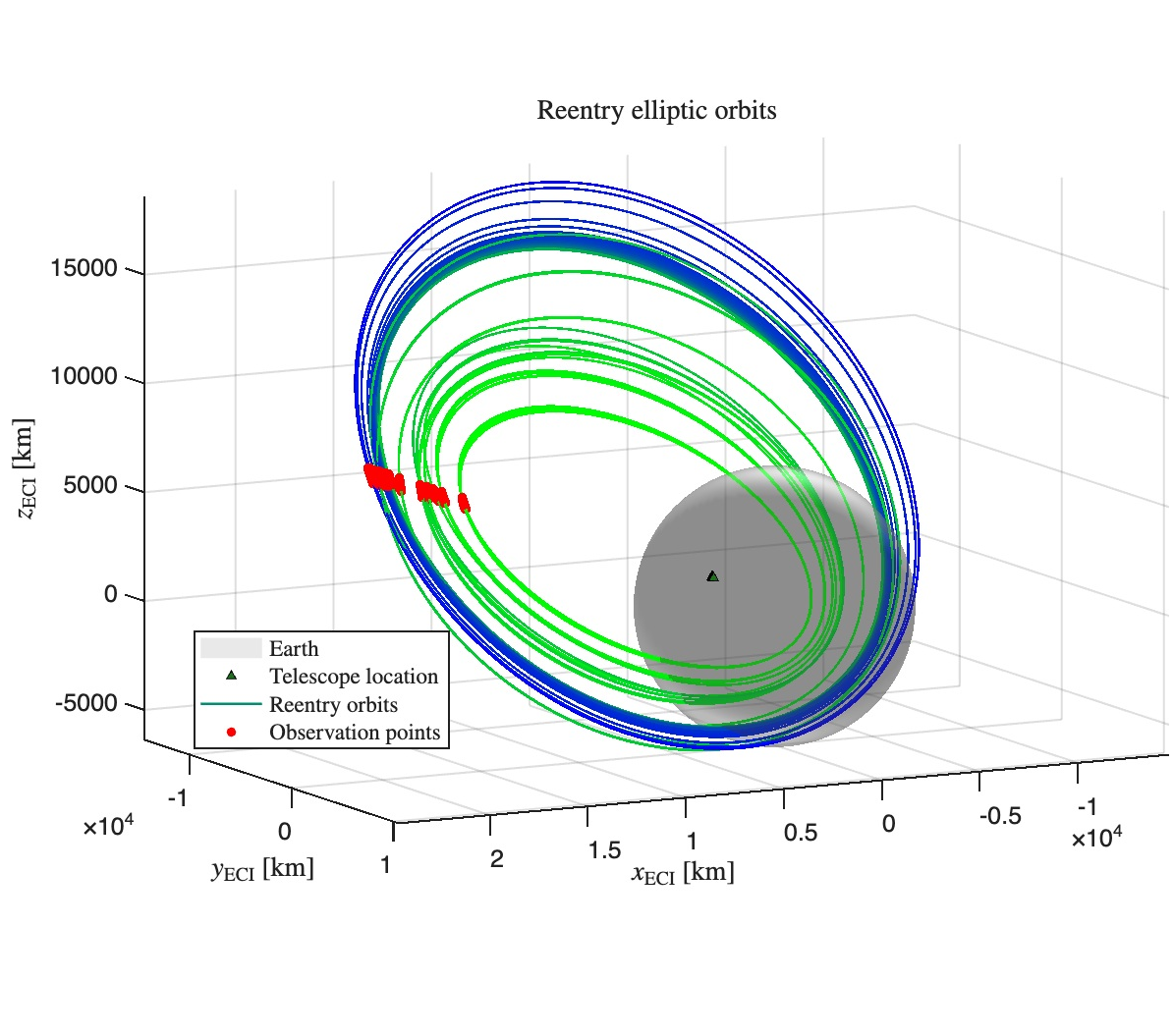}
    \label{fig:reentry}
}
\hfill
\subfigure[Bounded elliptic orbits.]{
    \includegraphics[width=0.48\textwidth]{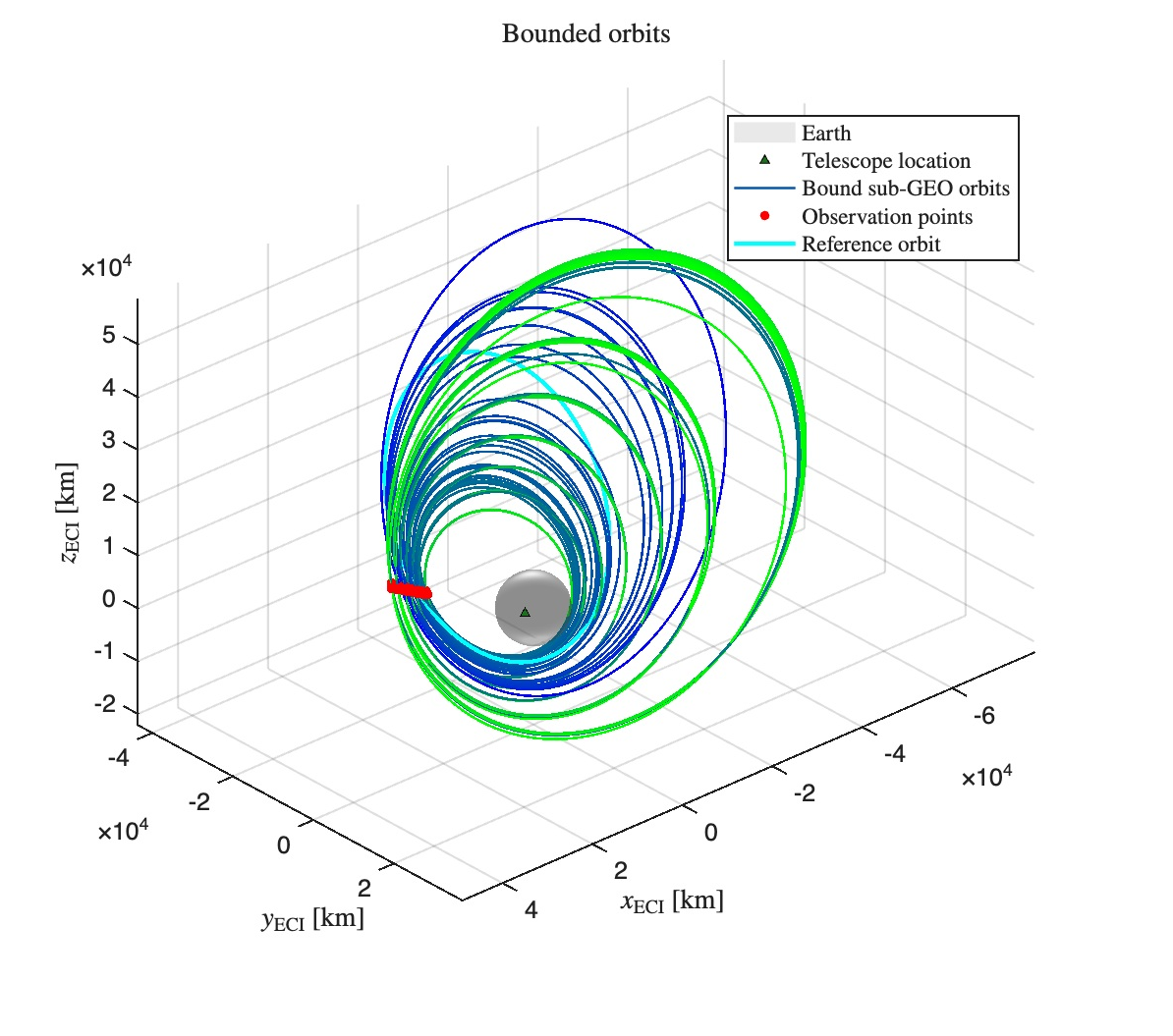}
    \label{fig:bounded}
}\\[1ex]

\subfigure[xGEO elliptic orbits.]{
    \includegraphics[width=0.48\textwidth]{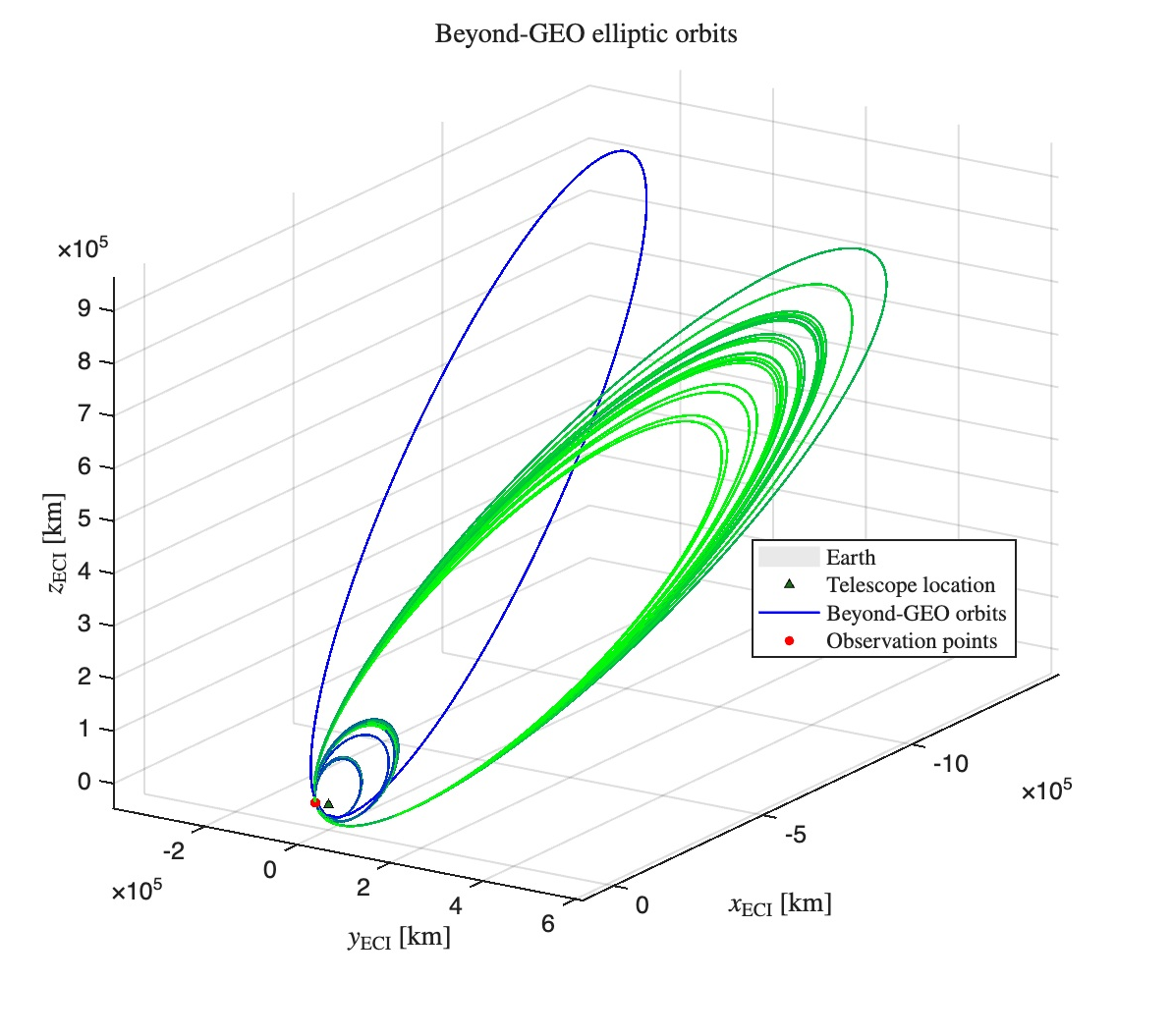}
    \label{fig:xGEO}
}
\hfill
\subfigure[Hyperbolic orbits.]{
    \includegraphics[width=0.48\textwidth]{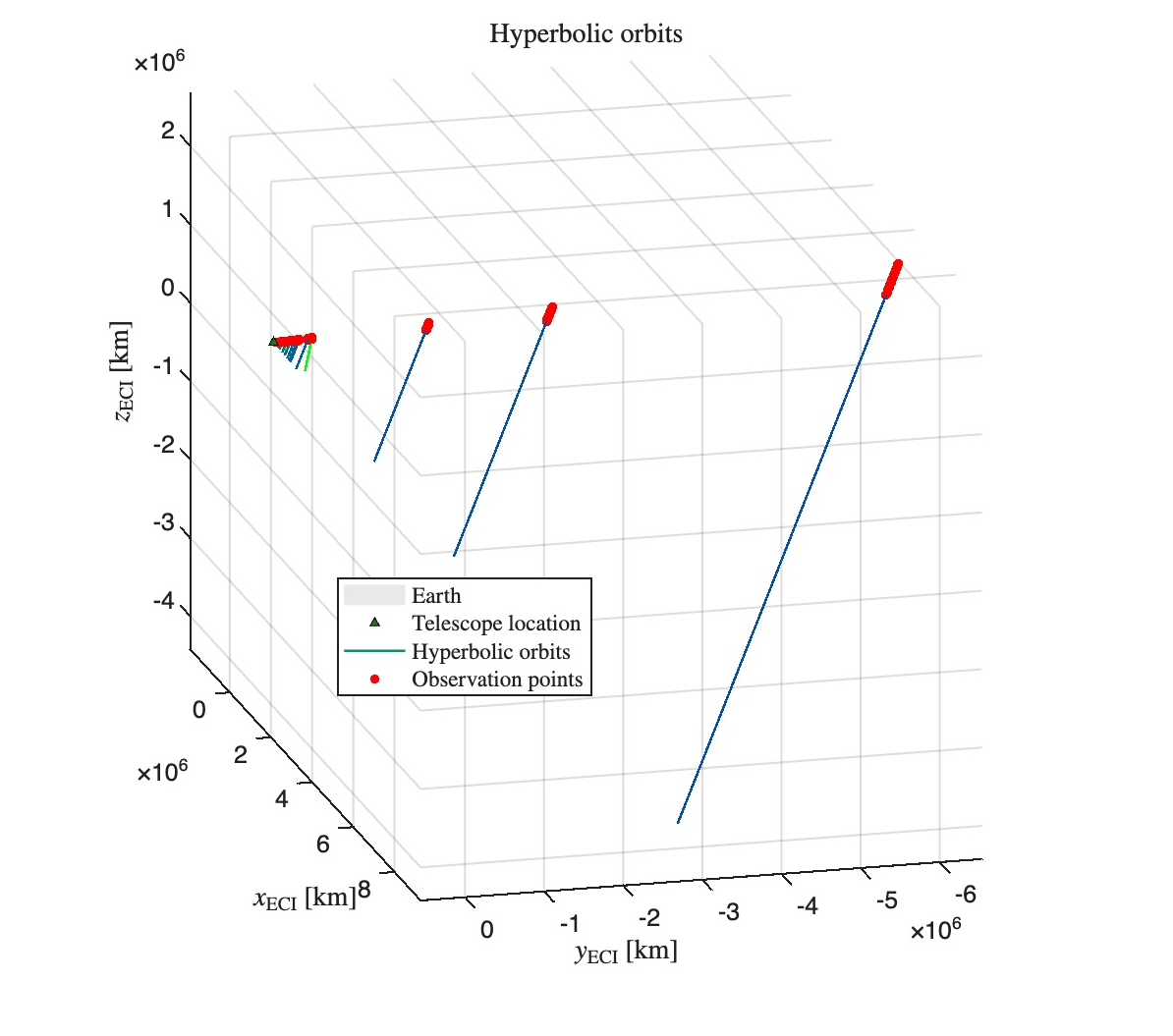}
    \label{fig:hyperbolic}
}

\caption{Three-dimensional view of the unconstrained solution classes. The candidate solutions are separated into reentry elliptic, bounded elliptic, xGEO elliptic, and hyperbolic orbit families. The color gradient denotes the final optimization cost, with blue corresponding to lower cost and green corresponding to higher cost. Red markers indicate the observation points. The reference Molniya orbit is shown in Figure~\ref{fig:bounded}.}
\label{fig:classes3d}

\end{figure}

The three-dimensional geometries of the unconstrained solution classes are shown in Figure~\ref{fig:classes3d}. The four subplots correspond to reentry elliptic orbits, bounded elliptic orbits, xGEO elliptic orbits, and hyperbolic orbits. The observation points are shown using red markers. The orbits are color-coded according to their final optimization cost, where blue corresponds to lower cost and green corresponds to higher cost.

A key observation from Figure~\ref{fig:classes3d} is that each class contains low-cost solutions. Thus, low residual cost alone does not uniquely identify the correct orbit family. This is a consequence of the short-arc angles-only geometry: different global orbits can pass through nearly the same line-of-sight directions during the observation interval.

\begin{figure}[htbp]
\centering
\subfigure[Close-Earth view.]{
    \includegraphics[width=0.48\textwidth]{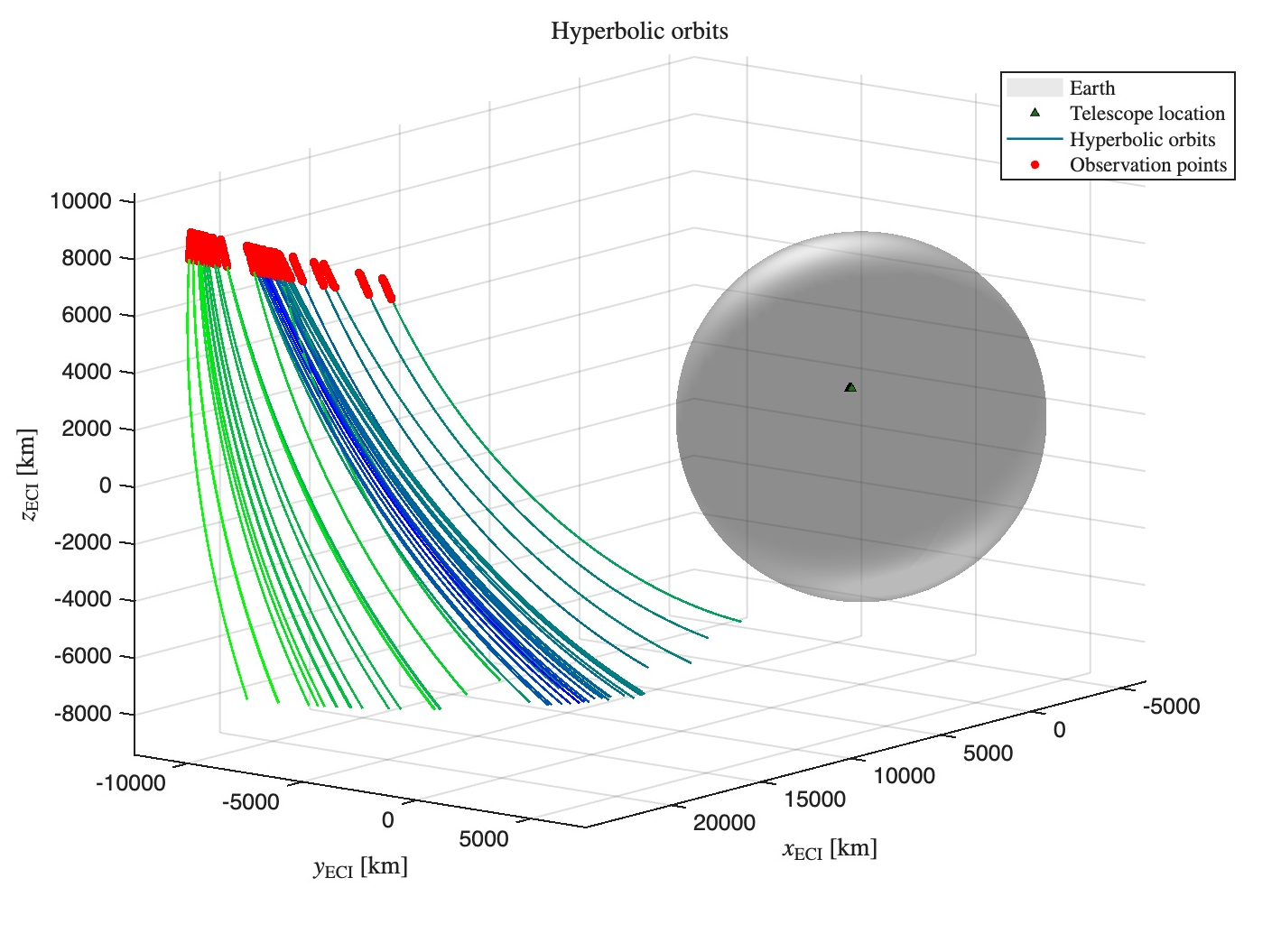}
    \label{fig:elliptic_unc_zoomed}
}
\hfill
\subfigure[Line-of-sight view.]{
    \includegraphics[width=0.48\textwidth]{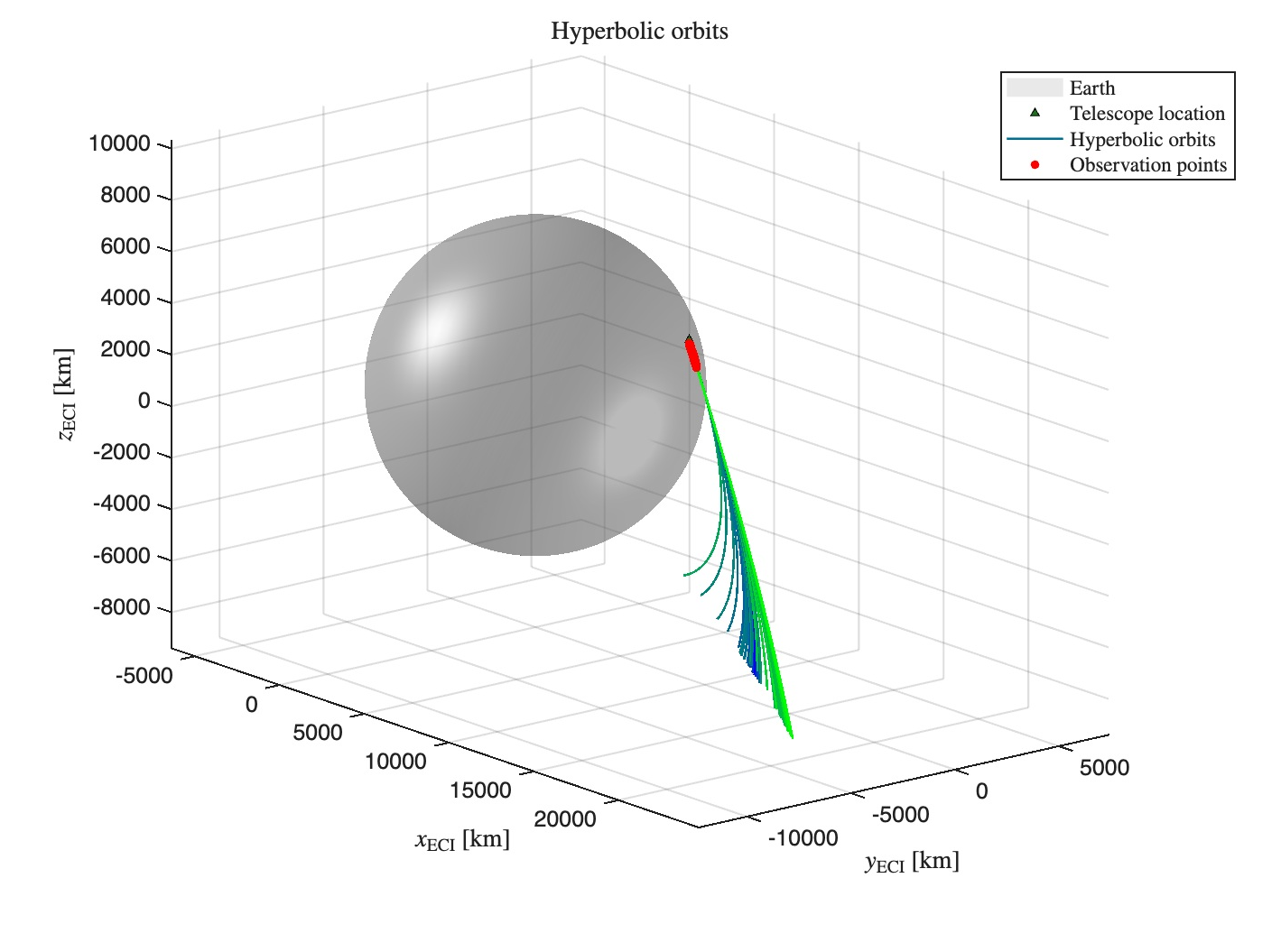}
    \label{fig:hyperbolic_unc_zoomed_los}
}

\caption{Close-Earth view of the hyperbolic solutions. The left subplot shows the local geometry of the open-orbit solutions near Earth, while the right subplot shows the corresponding line-of-sight view. The overlap of the trajectories near the observation directions explains why hyperbolic orbits can achieve low angular residuals over the short observation arc.}
\label{fig:hyperbolic_close_earth}

\end{figure}

The hyperbolic solutions are shown more clearly in Figure~\ref{fig:hyperbolic_close_earth}. Since the full hyperbolic trajectories extend over a large spatial region, their close-Earth behavior is difficult to interpret from Figure~\ref{fig:classes3d}. Figure~\ref{fig:hyperbolic_close_earth} therefore focuses on the region near Earth and near the observation arc. The line-of-sight view shows that the hyperbolic trajectories overlap closely near the observation directions. This explains why hyperbolic orbits can achieve low angular residuals over the short observation window, even though their global orbital behavior is very different from the Molniya reference orbit.

\begin{figure}[htbp]
\centering

\subfigure[All elliptic orbits.]{
    \includegraphics[width=0.48\textwidth]{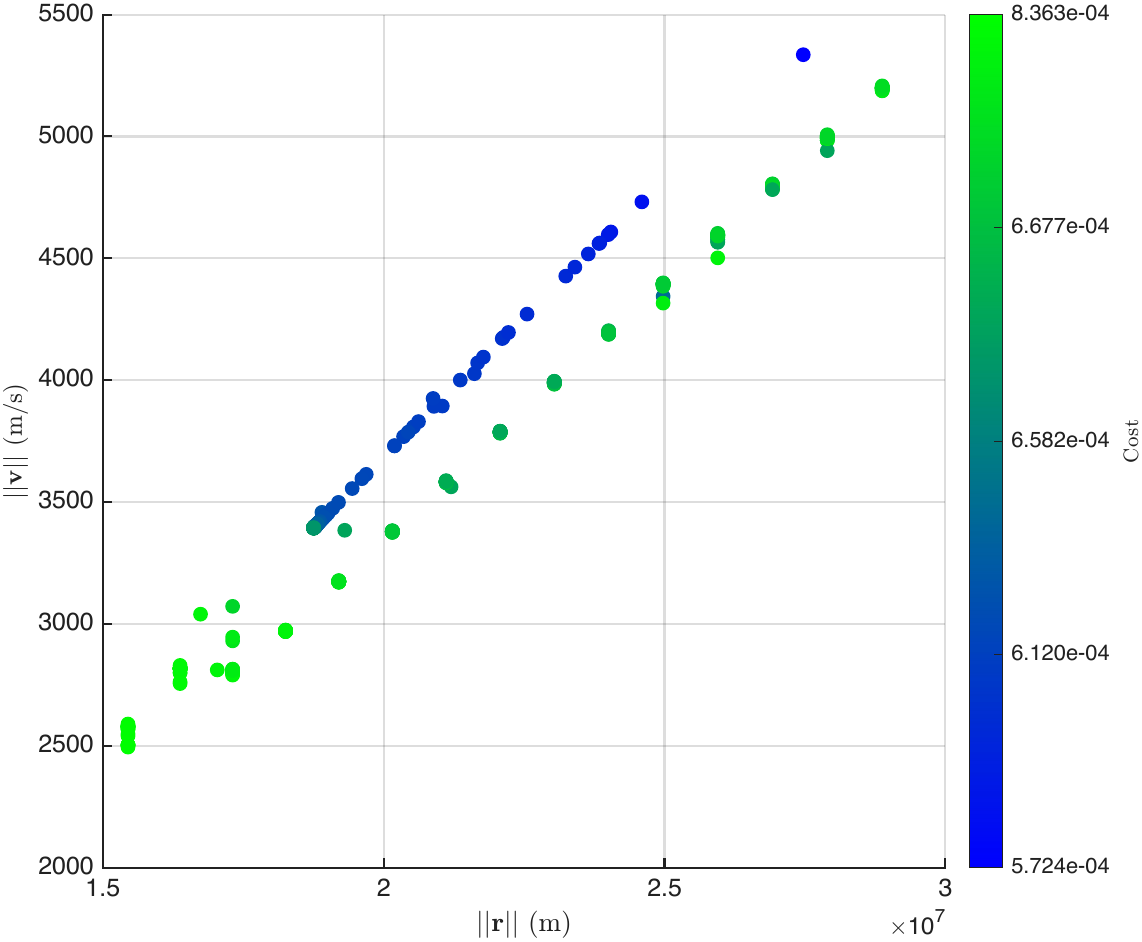}
    \label{fig:elliptic_unc_cost}
}
\hfill
\subfigure[All hyperbolic orbits.]{
    \includegraphics[width=0.48\textwidth]{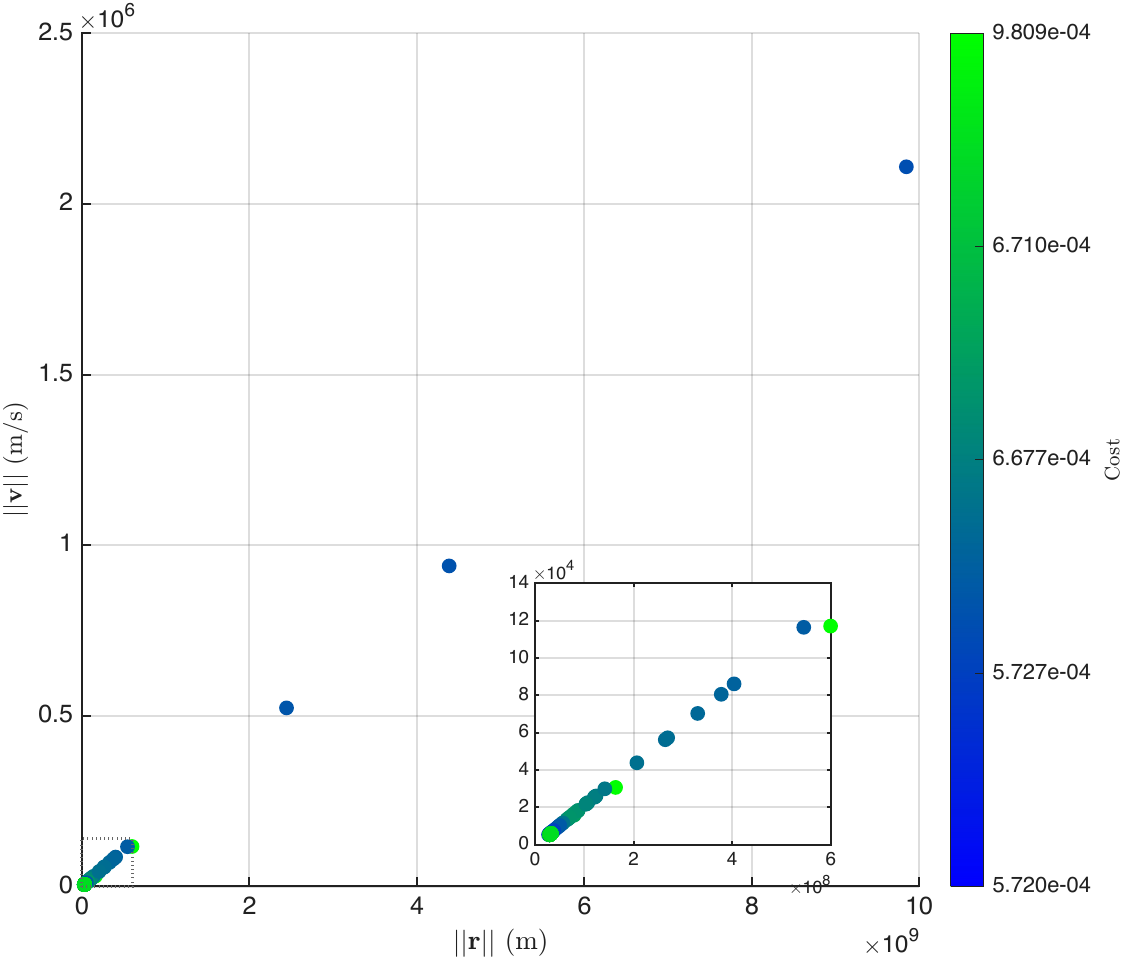}
    \label{fig:hyperbolic_unc_cost}
}

\caption{Cost-colored $||\mathbf{v}||$-versus-$(||\mathbf{r}||)$ representation of the unconstrained elliptic and hyperbolic solution families. Blue denotes lower-cost solutions and green denotes higher-cost solutions. The elliptic solutions separate into two visible cost branches, while the hyperbolic solutions retain an approximately linear structure with a zoomed inset showing the near-Earth region.}
\label{fig:elliptic_hyperbolic_cost}

\end{figure}

Figure~\ref{fig:elliptic_hyperbolic_cost} presents the same $(||\mathbf{v}||)$-versus-$(||\mathbf{r}||)$ classification as Figure~\ref{fig:elliptic_hyperbolic}, but with the final optimization cost included through a color scale. Blue denotes lower-cost solutions, while green denotes higher-cost solutions. The elliptic solutions separate into two visible cost branches: one branch contains lower-cost solutions, while the other contains higher-cost solutions. The hyperbolic plot also includes a zoomed inset to show the near-Earth region. The cost-colored view further confirms that the solution space contains multiple low-cost branches rather than a single isolated optimum.

\subsection{Constrained Optimization Results}

The unconstrained results demonstrate that the short-arc angles-only problem admits several low-cost orbit families. To isolate a physically plausible subset of solutions, the optimization problem was also solved using \texttt{fmincon} with the constraints described in the previous section. These constraints enforce elliptic motion, exclude reentry trajectories, and restrict the semimajor axis below the xGEO threshold.

It is important to emphasize that this constrained minimization is a heuristic. In a real angles-only short-arc orbit determination problem, one generally does not know a priori that the object must satisfy these constraints. Therefore, the constrained formulation should not be interpreted as proving uniqueness of the solution. Instead, it is used to examine the subset of solutions that are consistent with a desired physically admissible orbit class.

\begin{figure}[htbp]
\centering
\subfigure[Constrained orbit solutions.]{
    \includegraphics[width=0.48\textwidth]{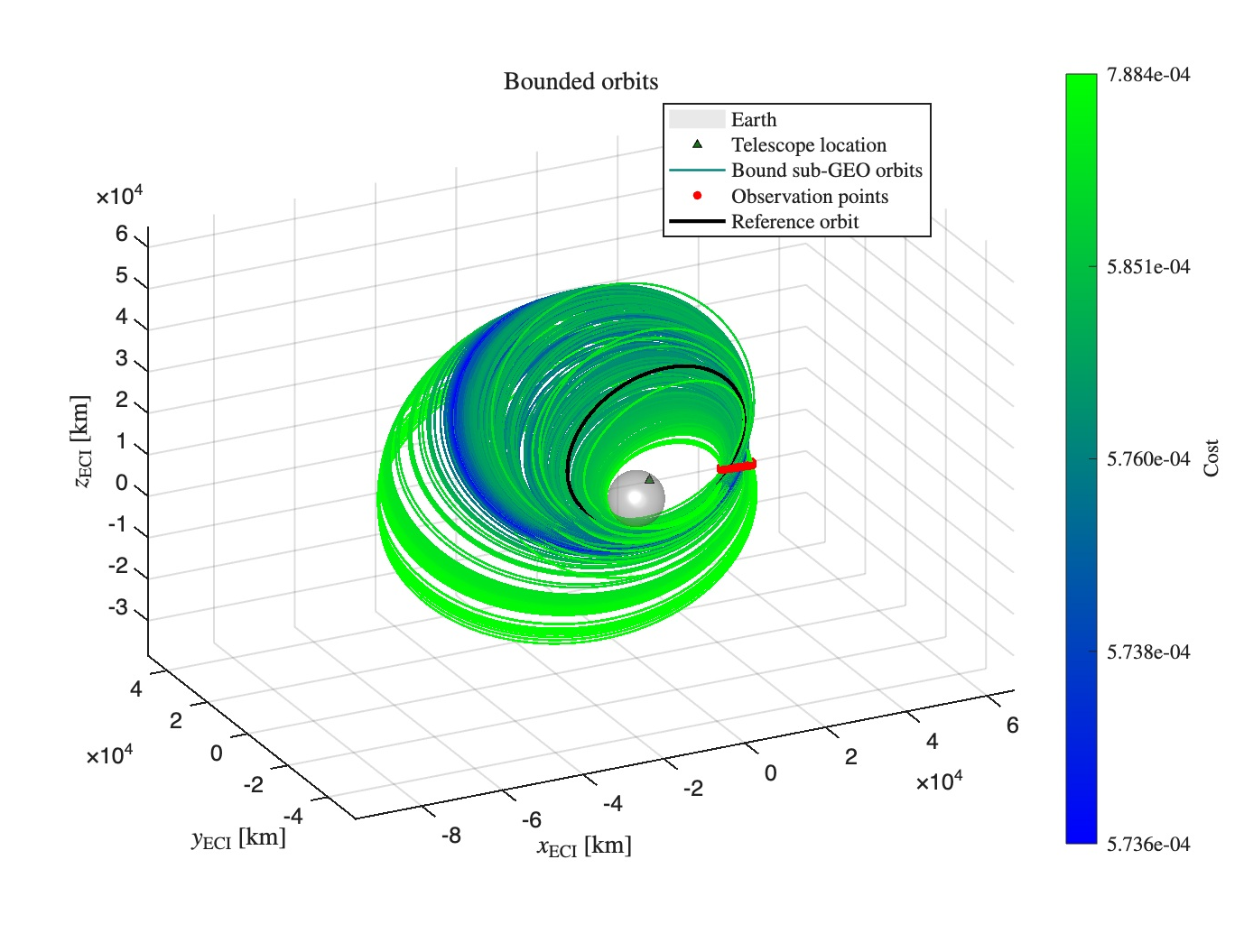}
    \label{fig:view_con}
}
\hfill
\subfigure[Line-of-sight view.]{
    \includegraphics[width=0.48\textwidth]{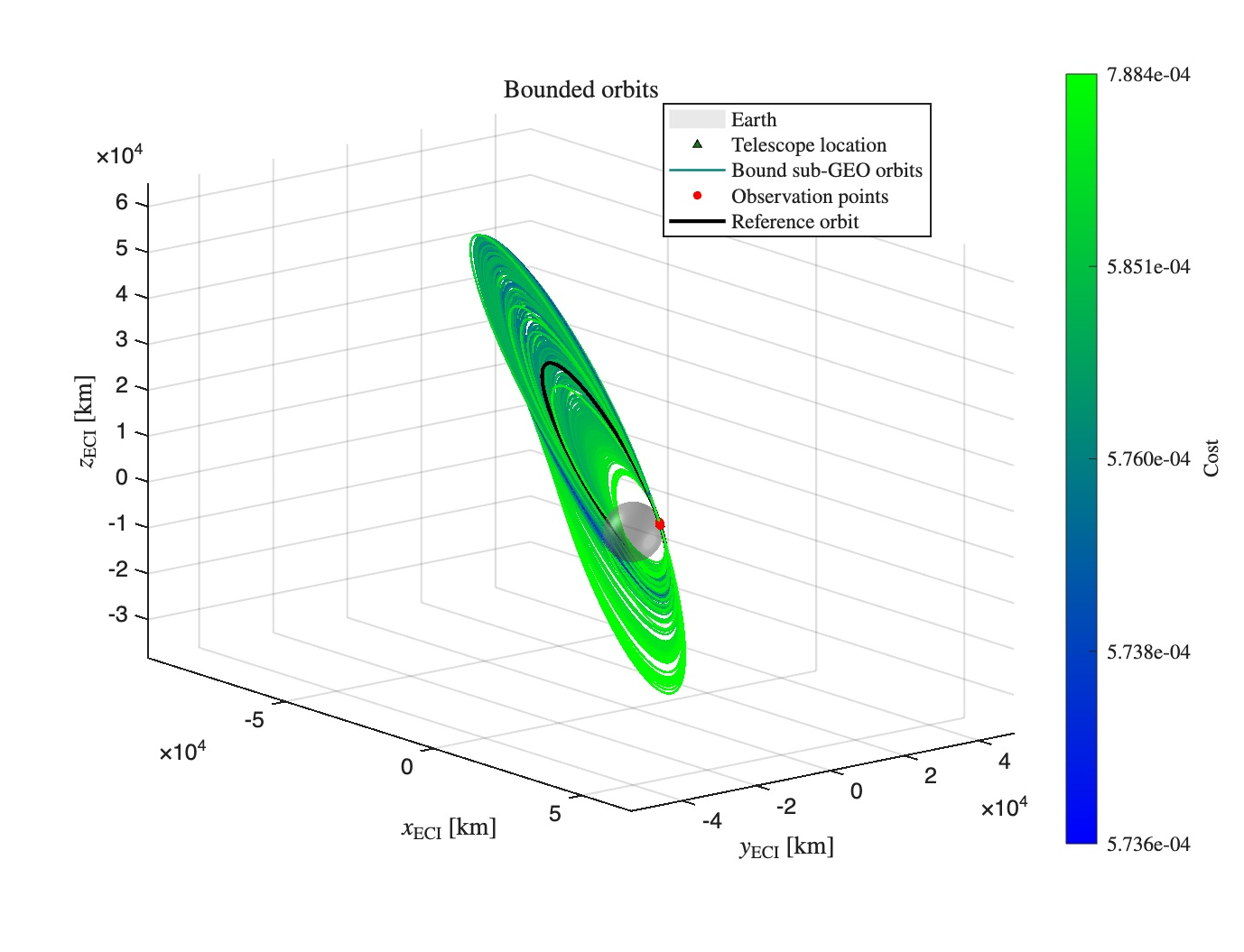}
    \label{fig:los_con}
}

\caption{Constrained \texttt{fmincon} solutions obtained after enforcing elliptic, non-reentry, and sub-xGEO constraints. The left subplot shows the resulting bounded elliptic trajectories, and the right subplot shows the line-of-sight view near the observation arc. The red markers denote the observation points. These constraints are used as a heuristic to isolate physically plausible Molniya-like solutions.}
\label{fig:elliptic_con}

\end{figure}

Figure~\ref{fig:elliptic_con} shows the solutions obtained from the constrained \texttt{fmincon} formulation. The first subplot shows the resulting bounded elliptic trajectories, while the second subplot shows the line-of-sight view near the observation arc. The red markers denote the observation points. The line-of-sight view shows that the constrained solutions remain closely aligned near the measurement directions, which again indicates that the short-arc angular data primarily constrain the local observation geometry. However, compared with the unconstrained case, the reentry, xGEO, and hyperbolic families are removed by construction.

The orbital-element range obtained from the constrained solutions is summarized in Table~\ref{tab:range_fmincon}. The table includes the true Molniya orbit parameters, the minimum and maximum values recovered from the constrained solution set, and a weighted solution. The weighted solution is computed using inverse-cost weights,
\[
w_j =
\frac{1/J_j}{\sum_k 1/J_k},
\]
where $(J_j)$ is the final cost of the (j)-th solution. Thus, lower-cost solutions contribute more strongly to the weighted estimate.

The weighted solution is close to the reference Molniya orbit. The inclination, right ascension of the ascending node, and the combined angle $(\omega+\nu_0)$ are especially close to the true values. As in the unconstrained case, the individual values of $(\omega)$ and $(\nu_0)$ show a wider variation, but their sum remains tightly constrained. This supports the conclusion that the short arc constrains the observed angular position more strongly than it constrains the individual orbital-angle decomposition.


\begin{table}[!htbp]
\centering
\begin{tabular}{|p{2.0cm}|p{2.2cm}|p{1.1cm}|p{1.3cm}|p{1.3cm}|p{1.3cm}|p{1.3cm}|p{1.4cm}|}
\hline
\textbf{Orbital Parameters} & a (m) & e & i & $\Omega$ & $\omega$ & $\nu_{0}$ & $\omega + \nu_{0}$ \\
\hline
\textbf{Actual Value} & $2.6559\times 10^{7}$ & 0.6757 & $63.5^{\circ}$ & $167.25^{\circ}$ & $277.98^{\circ}$ & $240^{\circ}$ & $517.98^{\circ}$\\
\hline
\textbf{Observed min.} & $1.3298\times 10^{7}$ & 0.21833 & $61.325^{\circ}$ & $158.69^{\circ}$ & $178.24^{\circ}$ & $185.41^{\circ}$ & $517.35^{\circ}$\\
\hline
\textbf{Observed max.} & $4.2\times 10^{7}$ & 0.71499 & $75.968^{\circ}$ & $169.1^{\circ}$ & $334.02^{\circ}$ & $342.2^{\circ}$ & $520.44^{\circ}$\\
\hline
\textbf{Weighted Sol.} & $2.4739\times 10^{7}$ & 0.50024 & $68.886^{\circ}$ & $163.82^{\circ}$ & $272.78^{\circ}$ & $246.35^{\circ}$ & $519.13^{\circ}$\\
\hline
\end{tabular}
\caption{Reference Molniya orbital parameters and the range of orbital elements recovered from the constrained \texttt{fmincon} solution set. The weighted solution is computed using inverse-cost weights and provides a representative estimate within the constrained family. The last two rows show the solutions obtained by using the approximate range estimate from the angular-rate analysis as the initial guess for the unconstrained and constrained optimizers, respectively.}
\label{tab:range_fmincon}
\end{table}

Overall, the results demonstrate that short-arc angles-only Molniya orbit determination is inherently ambiguous. The unconstrained optimization reveals reentry, bounded elliptic, xGEO elliptic, and hyperbolic solutions with comparable costs. The constrained optimization removes several undesirable solution families and produces a more Molniya-like set of solutions, but this improvement relies on heuristic constraints. Thus, the constrained formulation is useful for classification and for selecting physically plausible solutions, but it does not remove the fundamental multiple-solution nature of the short-arc angles-only orbit determination problem.

\section{Conclusions and Future Work}
\label{sec:conclusions}

This paper investigated the ambiguity of noisy short-arc angles-only orbit determination for a Molniya-orbit case study. The results show that classical Gauss IOD followed by GLSDC is unreliable under the considered short-arc and high-noise conditions, with several initial guesses failing to converge and the converged cases predominantly collapsing to hyperbolic solutions. To better explore the solution space, a Lambert-based initialization strategy was used over varying range assumptions and angular-observation pairs, followed by nonlinear optimization. The unconstrained results revealed that multiple distinct orbit families, including reentry elliptic, bounded elliptic, xGEO elliptic, and hyperbolic trajectories, can reproduce the same short-arc angular measurements with comparable residuals. This demonstrates that the final cost alone is insufficient to identify the correct orbit family. The recovered solutions also exhibited a structured range--velocity relationship, where larger estimated position magnitudes corresponded to larger velocity magnitudes, indicating that the angles-only data primarily constrain the apparent line-of-sight angular motion rather than the absolute range and range-rate. Physically motivated constraints were then introduced to isolate Molniya-like bounded elliptic solutions, and the inverse-cost weighted constrained solution was found to be close to the reference orbit. However, these constraints should be interpreted as heuristic a priori information rather than proof of uniqueness, since such information may not be available in a real short-arc orbit determination scenario.

Future work will focus on performing a larger Monte Carlo analysis over measurement noise, observation geometry, arc length, and initial orbital conditions in order to quantify the probability distribution over the candidate orbit families. In addition, instead of relying on \texttt{fmincon}, we will develop a customizable Newton-based optimization framework that allows more direct control over the Hessian approximation, line search, constraint handling, and convergence behavior. Future work will also incorporate additional information, such as apparent brightness or magnitude measurements, object-size and reflectivity assumptions, orbit-class constraints, and follow-up observations, to further restrict the solution space and provide a more restrictive characterization of the admissible candidate orbits.

\bibliographystyle{AAS_publication}   
\bibliography{references}   

\end{document}